# Disaggregation in Bundle Methods: Application to the Train Timetabling Problem

**Abderrahman Ait Ali** [a,1], Per Olov Lindberg [a,2], Jan-Eric Nilsson [3],
Jonas Eliasson [a,4], Martin Aronsson [5]
[a] Department of Transport Science, KTH Royal Institute of Technology
Teknikringen 10, 114 28 Stockholm, Sweden
[1] E-mail: abde@kth.se, Phone: +46 (0) 8 790 9439
[2] E-mail: polin@kth.se
[3] E-mail: jan-eric.nilsson@vti.se
[4] E-mail: Jonas.eliasson@abe.kth.se
[5] E-mail: martin@sics.se

**Abstract**
Bundle methods are often used to solve dual problems that arise from Lagrangian relaxations of large scale optimization problems. An example of such problems is the train timetabling problem. This paper focuses on solving a dual problem that arises from Lagrangian relaxation of a train timetabling optimization program. The dual problem is solved using bundle methods. We formulate and compare the performances of two different bundle methods: the aggregate method, which is a standard method, and a new, *disaggregate*, method which is proposed here. The two methods were tested on realistic train timetabling scenarios from the Iron Ore railway line. The numerical results show that the new disaggregate approach generally yields faster convergence than the standard aggregate approach.

## 1 Introduction

Several train timetabling models have been developed in various contexts, see e.g. (Caprara et al., 2006), (Jamili et al., 2012) and (Xu et al., 2014). Those models are often based on mathematical programming such as Integer Programming (IP) or Mixed Integer Programming (MIP) combined with solution methods that use the mathematical properties of the model.

One of the early and most cited papers to deal with this problem and research area is (Brännlund et al., 1998). The paper presents a MIP formulation of the train timetabling problem (TTP). The MIP model is NP-hard and thus difficult to solve for large applications of the problem. The authors use a Lagrangian relaxation solution approach where the track capacity constraints are relaxed and assigned prices or multipliers. A dual iterative scheme with a heuristic to find feasible solutions is used to solve the dual problem.

Several solution methods aggregate the different train requests while solving the dual problem. For instance, (Caprara et al., 2002) proposed a multigraph theoretical formulation of the problem to derive an integer linear programming model. The model was solved using Lagrangian relaxation of the constraints associated with the nodes in the graph. However, there was not any use of disaggregation between the different train requests.

It is possible to include many interesting functional constraints into the train timetabling model. (Caprara et al., 2006) developed a more realistic model by including several functionalities such as infrastructure maintenance, timetable prescription and manual block signalling. A Lagrangian heuristic algorithm was used to solve the problem without disaggregating the train requests.

(Schlechte, 2012) presents a group of models and algorithms for railway track allocation which includes integer programming models. The problem was formulated using the standard train scheduling graph. The model was later improved by taking advantage of the structure of the headway conflicts. The combinatorial problem was solved using a two-step solution method. The first step is to solve a linear relaxation of the problem using proximal bundle method without disaggregation. The second step is to solve the original problem using the linear relaxation solution with a sophisticated branch and price algorithm.

One property of the train timetabling problem which has yet not been thoroughly investigated, is the disaggregation property between the train requests. In this application to train timetabling, disaggregation uses a problem property in order to speed up the dual solution methods by dividing the aggregate problem into smaller sub-problems for each train request.

This paper starts from an aggregate formulation based on the model in (Brännlund et al., 1998) before deriving a novel bundle method approach to solve the dual problem, which arises from the relaxation of (TTP). This novel approach uses the disaggregation property of the problem. We compare this *disaggregate bundle method* with the traditional aggregate approach. We first present the problem, the notations and the mathematical model. Second, we derive and describe the solution approach where we distinguish between the aggregate and disaggregate approach. Finally, the comparative numerical results are given and relevant conclusions are drawn.

## 2 Mathematical Model

### 2.1 Problem Description

The network that is considered in the train timetabling problem is a single track line. This does not affect the generality of the problem. The line is discretized into different blocks. The blocks are of two types: *station blocks* where the trains are allowed to meet and *signalling blocks* where they are not. Each block of the line has a certain capacity that represents the number of tracks in that particular block. In particular, station blocks have a capacity higher or equal to one whereas signalling blocks have a capacity of one. Figure 1 illustrates the two types of blocks in a stretch of a single-track line.

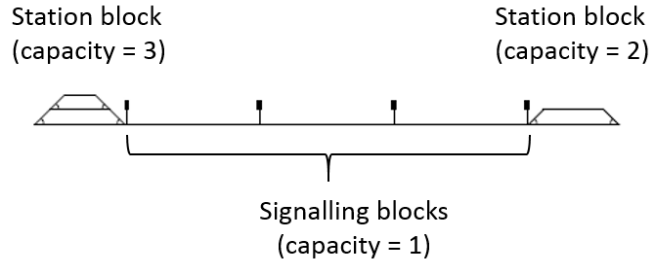

Figure 1- A stretch of a single-track rail line

In a station block, the train can either pass at full speed if no compulsory stop is requested, or stop before leaving for the next destination. In a signalling block, the train can pass without stopping or wait. This leads to four different speed scenarios in each block depending on the state of the train at the entrance or exit of the block: FF, SF, SS, FS where F is full speed state and S is stopping state. In Figure 2, the four scenarios are illustrated; FF is the fastest scenario and SS is the slowest. SF and FS are slower than FF but faster than SS.

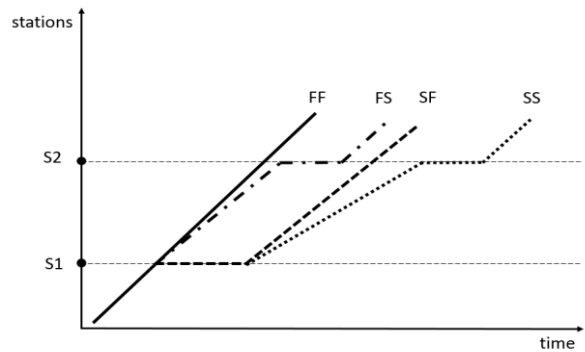

Figure 2 -The four different train speed scenarios between two train stations

In these different scenarios, the acceleration and deceleration (i.e. braking) properties of the trains play an important role in determining the travel time for the different scenarios.

## 2.2 Notations

The time is discretized into *time intervals* denoted $t \in \mathcal{T} = \{1, \dots, T\}$. Similarly, the single track line is discretized into *blocks* labelled $b \in \mathcal{B} = \{1, \dots, B\}$. The *train requests* for traffic are denoted $r \in \mathcal{R} = \{1, \dots, R\}$.

It is possible to see this problem in a graph perspective. Indeed, for each train request r, the different possibilities to perform its traffic is given by a network $\mathcal{N}_r = (\mathcal{A}_r, \mathcal{V}_r)$, where nodes $\mathcal{V}_r$ represent positions in space and time ($\mathcal{B} \times \mathcal{T}$), and whose arcs $\mathcal{A}_r$ represent train movements. The same graph perspective is also mentioned in (Brännlund et al., 1998). For instance, one arc could be "leaving block $b_1$ at time t, accelerating from standstill, travelling to block $b_2$ without retarding there", or "standing still at block $b_1$ from time $t$ to $t + 1$".

Each train request r thus has a finite set of feasible paths $p \in \mathcal{P}_r$ (from a starting station $S_r^s$ through the network to an ending station $S_r^e$) to perform its duties. A path p is an ordered subset of $\mathcal{A}_r$ that describes the trajectory of the train in time and space. The set of all train paths is $\mathcal{P} = \cup_r \mathcal{P}_r$. Furthermore, we assume that one of the feasible paths is the "null path" that corresponds to not scheduling the train request.

For a certain path $p \in \mathcal{P}_r$ for a train request $r \in \mathcal{R}$, the train movement of arc $a \in p \subseteq \mathcal{A}_r$ leads to the (not necessarily physical) occupation of a certain block-times, given by the matrix $\Delta^a = (\delta_{bt}^a)$ with

$$\delta_{b,t}^a = \begin{cases} 1 & \text{if movement } a \text{ occupies blocktime } (b,t) \\ 0 & \text{otherwise} \end{cases}. \quad (1)$$

To each path $p \in \mathcal{P}_r$ is associated a utility value $v_p$ of choosing the path. We will return to the properties of this value in Section 4.2.

The following summarizes the different notations that will be used in this paper.

- Sets

| | |
|---|---|
| $\mathcal{T}$ | Set of time intervals $\{1, 2, \ldots, t, \ldots, T\}$. |
| $\mathcal{B}$ | Set of space blocks $\{1, 2, \ldots, b, \ldots, B\}$. |
| $\mathcal{R}$ | Set of train requests $\{1, 2, \ldots, r, \ldots, R\}$. |
| $\mathcal{P}_r$ | Set of paths for train request $r \in \mathcal{R}$. |
| $\mathcal{P} = \cup_r P_r$ | Set of all train paths. |

- Parameters

| | |
|---|---|
| $c_b$ | Capacity of space block $b$. |
| $\delta_{b,t}^a$ | Capacity usage in time $t$ of block $b$ for train movement $a$. |
| $v_p$ | Value of the allocation of path $p$. |

- Decision variables

| | |
|---|---|
| $x_p \in \{0,1\}$ | Allocation state of path $p \in \mathcal{P}$. |

## 2.3 Model

The TTP is now to select one path p for each train request r, such that capacity limits are not violated, and such that the total timetable utility is maximized. By introducing a binary indicator variable $x_p$ for the selection/non-selection of path p to be run, the problem can be stated as follows:

$$(TTP) \begin{cases} \max_{x_p} \sum_{p \in \mathcal{P}} v_p x_p \\ s.t. \begin{cases} \sum_{p \in \mathcal{P}} (\sum_{a \in p} \delta_{b,t}^a) x_p \leq c_b, & \forall (b,t) \in \mathcal{B} \times \mathcal{T} \quad (i) \\ \sum_{p \in \mathcal{P}_r} x_p = 1, & \forall r \in \mathcal{R} \quad (ii) \\ x_p \in \{0,1\}, & \forall p \in \mathcal{P} \quad (iii) \end{cases} \end{cases} \quad (2)$$

Constraints (2.i) are the capacity constraints of the blocks, (2.ii) the constraints to choose exactly one path for each train request, and (2.iii) are the binary constraints on the path selection variable $x_p$.

The factor $(\sum_{a \in p} \delta_{b,t}^a)$ in front of $x_p$ in (2.i) corresponds to the capacity consumption of path p on block-time (b, t). It can be denoted as $d_{bt}^p$ and therefore regrouped in the matrix $d^p = (d_{bt}^p)$. Therefore, constraints (2.i) become

$$\sum_{p \in \mathcal{P}} d_{bt}^p x_p \leq c_b, \quad \forall (b,t) \in \mathcal{B} \times \mathcal{T}. \tag{2.i'}$$

We will use this new form in all the formulations that follow.

## 3 Solution Methods

(TTP) is an IP with a very large number of binary variables. The combinatorial nature of the problem makes it extremely difficult to solve for large instances using the-state-of-the-art IP solvers. In ordered to get around the computational complexity of the problem, we use the classical Lagrangian relaxation technique as the starting point for the solution method.

In the Lagrangian relaxation of (TTP), we are allowed to violate constraints (2.i'), but at a certain price given by the corresponding Lagrangian multipliers $\mu = (\mu_{bt}) \geq 0$. The relaxed (TTP), noted $(TTP)_\mu$, is formulated as

$$(TTP)_\mu \begin{cases} \varphi(\mu) := \max \sum_{p \in \mathcal{P}} v_p x_p + \sum_{(b,t) \in \mathcal{B} \times \mathcal{T}} \mu_{bt}(c_b - \sum_{p \in \mathcal{P}} d_{bt}^p x_p) \\ s.t. \begin{cases} \sum_{p \in \mathcal{P}_r} x_p = 1, & \forall r \in \mathcal{R} \\ x_p \in \{0,1\}, & \forall p \in \mathcal{P} \end{cases} \end{cases} \tag{3}$$

$(TTP)_\mu$ is a relaxation of (TTP) since every feasible solution to (TTP) is also a feasible solution to $(TTP)_\mu$, and the objective value of any feasible solution in (TTP) is not greater than that in $(TTP)_\mu$. Hence, for each value of $\mu$, the objective value $\varphi(\mu)$ in $(TTP)_\mu$ is larger than or equal to the optimal value of (TTP).

It is possible to further simplify the objective value of $(TTP)_\mu$. For a given $\mu \geq 0$ and under the same constraints as in (3), $\varphi(\mu)$ can be rewritten as

$$\varphi(\mu) = \sum_{(b,t) \in \mathcal{B} \times \mathcal{T}} c_b \mu_{bt} + \max_{x_p} \sum_{p \in \mathcal{P}} (v_p - \sum_{(b,t) \in \mathcal{B} \times \mathcal{T}} \mu_{bt} d_{bt}^p) x_p. \tag{4}$$

The factor in front of $x_p$ in (4) or the term $(v_p - \sum_{(b,t) \in \mathcal{B} \times \mathcal{T}} \mu_{bt} d_{bt}^p)$, can be considered as a reduced utility revenue from choosing path $p \in \mathcal{P}$ (i.e. $x_p := 1$) under the current multipliers $\mu$. Thus, $(TTP)_\mu$ decomposes into one shortest-path problem for each train request $r \in \mathcal{R}$. Such shortest path problems have well-established solution algorithms and are relatively easy to solve. In this model, we developed a shortest path algorithm based on topological sorting. This is justified by the fact that the train movements network is a weighted directed acyclic graph (Cormen et al., 2009).

### 3.1 Dual Problem

As just noted, $\varphi(\mu)$ is larger than or equal to the optimal value of the original problem (TTP), for any $\mu \geq 0$. It is therefore an upper bound to the optimal value of (TTP). Thus, the dual problem (D) is to find the optimal solution $\mu^*$ that gives the best (i.e. smallest) upper bound.

$$\text{(D)} \quad \begin{cases} \min \varphi(\mu) \\ s.t. \quad \mu \geq 0 \end{cases} \tag{5}$$

Since there are only a finite number of shortest path combinations, $\varphi$ is piecewise linear (or rather piecewise affine). It is therefore a convex function since it is the maximum of a set of linear functions. Moreover, $\varphi$ has a lower bound, i.e. any feasible solution to the original problem (TTP). Therefore, (D) has a global minimum $\varphi^*$ at the optimal multipliers $\mu^*$.

Let us assume that for an arbitrary value $\bar{\mu} \geq 0$, the maximum in $(TTP)_{\bar{\mu}}$ is achieved at $\tilde{x}(\bar{\mu}) = (\tilde{x}_p)_{p \in \mathcal{P}}$. Inserting the corresponding $\tilde{x}(\bar{\mu})$ in the objective of $(TTP)_\mu$ gives a linear (in $\mu$) function $\tilde{\varphi}(\mu) = \sum_{p \in \mathcal{P}} v_p \tilde{x}_p + \sum_{bt} \mu_{bt}(c_b - \sum_{p \in \mathcal{P}} d_{bt}^p \tilde{x}_p)$, that is equal to $\varphi(\mu)$ at $\bar{\mu}$. This linear function corresponds to a supporting plane to the graph of $\varphi$. The slope of this function is given by the matrix $g(\bar{\mu}) = (\bar{g}_{bt}) \coloneqq (c_b - \sum_{p \in \mathcal{P}} d_{bt}^p \tilde{x}_p)$ which is a subgradient of $\varphi$ at $\bar{\mu}$. Thus the supporting linear function to $\varphi$ at $\bar{\mu}$ can be written as

$$\varphi(\bar{\mu}) + g(\bar{\mu}) * (\mu - \bar{\mu}), \tag{6}$$

where * denotes the inner product between two matrices (i.e. component wise).

In order to solve (D), we use the aggregate bundle method described in (Kiwiel, 1990). Based on this method, a novel disaggregate approach is derived that uses the disaggregation property in the train timetabling problems.

### 3.2 The Aggregate Bundle Method

For $\mu = \bar{\mu}$, the (possibly many) $\tilde{x}(\bar{\mu})$ giving the maximum in $(TTP)_\mu$ give the subgradients to $\varphi$ at $\bar{\mu}$. Suppose that we currently are at $\mu = \mu_k$, and that we have chosen to approximate $\varphi$ by the supporting planes computed in iterations $l \in \mathcal{L}_k$, where $\mathcal{L}_k$ is the *bundle* at iteration $k$ from the previous iterations. Let the corresponding subgradients be $\{g_l\}_{l \in \mathcal{L}_k}$. Then in the standard aggregate bundle method, we compute a new tentative solution as the solution to the following sub-problem

$$(\overline{D}_k^{agg}) \quad \begin{cases} \min \bar{\varphi}^k(\mu) + \dfrac{u_k}{2} |\mu - \mu_k|^2 \\ s.t. \quad \mu \geq 0, \end{cases} \tag{7}$$

where $|\cdot|$ denotes the Euclidean norm (or 2-norm) of a matrix reshaped into a vector, and $\bar{\varphi}^k(\mu) \coloneqq \max_{l \in \mathcal{L}_k} \{\varphi(\mu_l) + g_l * (\mu - \mu_l)\}$ is the maximum of the supporting linear functions at $\mu_l$, for $l \in \mathcal{L}_k$, giving an outer linearization of $\varphi$. The quadratic second term helps avoid taking too large steps and the step size is adjusted using the control parameter $u_k$ at each iteration.

In order to get around the inner maximization, $(\overline{D}_k^{agg})$ can be formulated as a single minimization problem by adding an additional variable as well as new constraints for the supporting linear functions. This leads to the following equivalent problem

$$(\overline{D}_k^{agg}) \begin{cases} \min v + \dfrac{u_k}{2}|\mu - \mu_k|^2 \\ s.t. \begin{cases} v \geq \varphi(\mu_l) + g_l * (\mu - \mu_l), & \forall l \in \mathcal{L}_k \quad (i) \\ \mu \geq 0 & (ii) \end{cases} \end{cases} \tag{8}$$

The matrices $\boldsymbol{\mu}_l$ for $l \in \mathcal{L}_k$ can be extremely large and lead to an excessive memory usage. Therefore, we suggest an equivalent formulation of the supporting linear functions in which scalars are stored instead of the matrices. We define this scalar, at $k$ and for all $l \in \mathcal{L}_k$, as

$$\Psi_{kl} := \varphi(\mu_l) + g_l * (\mu_k - \mu_l). \tag{9}$$

Hence, the problem in (7) becomes

$$(\overline{D}_k^{agg}) \begin{cases} \min v + \dfrac{u_k}{2}|\mu - \mu_k|^2 \\ s.t. \begin{cases} v \geq \Psi_{kl} + g_l * (\mu - \mu_k), & \forall l \in \mathcal{L}_k \quad (i) \\ \mu \geq 0 & (ii) \end{cases} \end{cases} \tag{10}$$

The $\Psi_{k+1,l}$ has to be updated recursively whenever the multipliers are updated, which means when $\mu_{k+1} \neq \mu_k$. The update scheme is

$$\Psi_{k+1.l} = \Psi_{k.l} + g_l * (\mu_{k+1} - \mu_k), \tag{11}$$

so that the supporting linear functions always have the current $\boldsymbol{\mu}_k$ as "foot point".

Let $y_{k+1}$ be the optimal solution to $(\overline{D}_k^{agg})$. At $y_{k+1}$ we evaluate the dual objective $\varphi$ by solving $(TTP)_\mu$ for $\mu = y_{k+1}$. We might then get a new supporting plane, including a new subgradient $g_{k+1}$. We define the achieved descent as $\bar{\varphi}^k(\mu_k) - \varphi(y_{k+1})$ and the forecasted one as $\bar{\varphi}^k(\mu_k) - \bar{\varphi}^k(y_{k+1})$. If the ratio of the achieved descent by the forecasted one is larger than a certain step quality threshold $m_L \in (0,1)$ then we set $\mu_{k+1} = y_{k+1}$, and the new $\mathcal{L}_{k+1}$ will incorporate the active supporting planes from $\mathcal{L}_k$ as well as the newly generated supporting plane. If otherwise the ratio is not large enough, we set $\mu_{k+1} = \mu_k$ and $\mathcal{L}_{k+1}$ will only add the newly generated supporting plane to $\mathcal{L}_k$. Thus, the polyhedral approximation of $\varphi$ is improved at each iteration.

The step control parameter $u_{k+1}$ is adjusted in both cases. It is set so that the curvature of the objective in $(\overline{D}_k^{agg})$ between $\mu_k$, and $y_{k+1}$ fits that of $\varphi$. The parameter has a minimum value, and is never decreased by a factor of more than 10 as in (Kiwiel, 1990)

### 3.3 The Disaggregate Bundle Method

In the disaggregated bundle method, we squeeze more information out of the solutions to the subproblems $(TTP)_\mu$ by considering a disaggregate dual objective. The dual objective function can be separated into independent functions for each train request.

$$\varphi(\mu) = \sum_{r \in \mathcal{R}} \varphi_r(\mu) + \sum_{(b,t) \in \mathcal{B} \times \mathcal{T}} c_b \mu_{bt}, \tag{12}$$

where $\varphi_r$ is defined for each train request $r \in \mathcal{R}$ as

$$\varphi_r(\mu) := \max \sum_{p \in \mathcal{P}_r} (v_p - \sum_{(b,t) \in \mathcal{B} \times \mathcal{T}} \mu_{bt} d_{bt}^p) x_p$$

$$s.t. \begin{cases} \sum_{p \in \mathcal{P}_r} x_p = 1, & (i) \\ x_p \in \{0,1\}, & \forall p \in \mathcal{P}_r \quad (ii) \end{cases} \tag{13}$$

which is the maximum income of train request $r \in \mathcal{R}$ under the current multipliers $\mu$. The multipliers can be interpreted as train access charges of using the rail infrastructure, an interpretation that we will return to in the concluding section.

The disaggregate dual problem, noted $(D^{dis})$, will be slightly different from (D) that is used in the aggregate approach.

$$(D^{dis}) \begin{cases} \min \sum_{r \in \mathcal{R}} \varphi_r(\mu) + \sum_{(b,t) \in \mathcal{B} \times \mathcal{T}} \mu_{bt} c_b \\ s.t. \quad \mu \geq 0 \end{cases}. \tag{14}$$

At iteration k in the disaggregate approach, each objective component $\varphi_r$, has its own bundle $\mathcal{L}_k^r$, with the subgradients defined as $g_{rl} := -d^{\hat{p}_{rl}}$, where $\hat{p}_{rl} \in \mathcal{P}_r$ is the shortest path in the sense that it leads to the maximal revenue for $\mu = \mu_l$ in $(TTP)_\mu$.

As with the aggregate approach, we use the subgradients to build supporting linear functions that are used as an outer approximation. In the disaggregate approach the outer approximation is computed for each objective component $\varphi_r$. Thus, the disaggregate bundle method problem is written as follows:

$$(\overline{D}_k^{dis}) \begin{cases} \min \sum_{r \in \mathcal{R}} v_r + \frac{u_k}{2} |\mu - \mu_k|^2 \\ s.t. \begin{cases} v_r \geq \varphi_r(\mu_l) + g_{rl} * (\mu - \mu_l), & \forall l \in \mathcal{L}_k \, \forall r \in \mathcal{R} \quad (i) \\ \mu \geq 0 & (ii) \end{cases} \end{cases}. \tag{15}$$

In order to minimize the memory storage in the implementation, as with the aggregate approach, instead of storing all the previous matrices of multipliers $\mu_l$, we only store the corresponding scalar parameters $\Psi_{r,l}^k$ that we define similarly as follows

$$\Psi_{r,l}^k := \varphi^r(\mu_l) + g_{rl} * (\mu_k - \mu_l). \tag{16}$$

The parameters are updated in a similar way to the aggregate case that is previously described. Thus, the formulation in (14) can be rewritten using the scalar parameters. So, at iteration $k$, we have

$$(\overline{D}_k^{dis}) \begin{cases} \min \sum_{r \in \mathcal{R}} v_r + \sum_{(b,t) \in \mathcal{B} \times \mathcal{T}} c_b \mu_{bt} + \frac{u_k}{2} |\boldsymbol{\mu} - \boldsymbol{\mu_k}|^2 \\ \text{s.t.} \begin{cases} v_r \geq \Psi_{r,l}^k + \mathbf{g_{rl}} * (\boldsymbol{\mu} - \boldsymbol{\mu_k}), \forall\, l \in \mathcal{L}_k^r, \forall\, r \in \mathcal{R} & (i). \\ \boldsymbol{\mu} \geq 0 & (ii) \end{cases} \end{cases} \quad (17)$$

We can further simplify (16) by introducing $\boldsymbol{s} = \boldsymbol{\mu} - \boldsymbol{\mu_k}$ to get the following formulation

$$(\overline{D}_k^{dis}) \begin{cases} \min \sum_{r \in \mathcal{R}} v_r + \sum_{(b,t) \in \mathcal{B} \times \mathcal{T}} (\mu_{bt}^k + s_{bt}) c_b + \frac{u_k}{2} |\mathbf{s}|^2 \\ \text{s.t.} \begin{cases} v_r \geq \Psi_{r,l}^k + \mathbf{g_{rl}} * \mathbf{s}, & \forall\, l \in \mathcal{L}_k^r, \forall\, r \in \mathcal{R} & (i) \\ \mathbf{s} \geq -\boldsymbol{\mu^k} & (ii) \end{cases} \end{cases} \quad (18)$$

## 4 Experimental Setup and Results

### 4.1 Implementation

Both solution methods, i.e. the aggregate and the disaggregate bundle method, are developed in MATLAB. The methods call a C++ program that computes the shortest path given the prices $\mu$ of occupying the block-times. The information between the two programming environments is exchanged using *mex functions* which are subroutines that allow MATLAB programs to call C, C++ and Fortran programs (MathWorks, 2016). In order to speed up the computation of the shortest path algorithm, the paths networks or the graphs of possible train movements are constructed once, and are stored in the C++ environment memory for use in all the iterations of the bundle method. Therefore, the MATLAB program (i.e. bundle method) calls at first a C++ function that allocates memory and constructs the train movement graphs from the input data before performing the bundle iteration. Figure 3 gives an overview of the software architecture that was implemented.

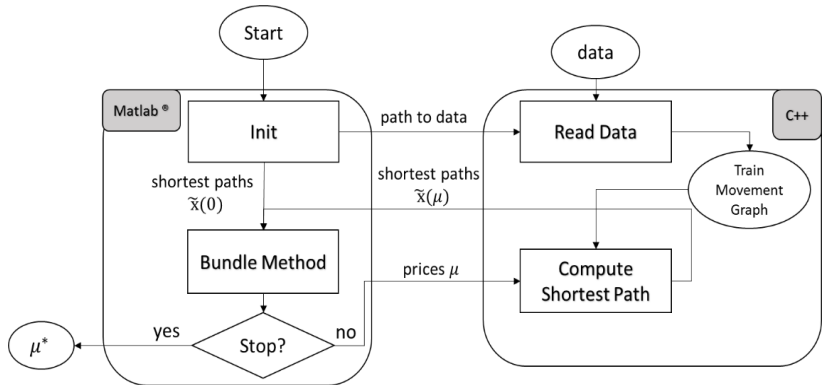

Figure 3 - Software architecture of the model implementation

## 4.2 Input Data

The input data that is used to test the two implementations is based on the information from train operations of the Iron Ore line (*Malmbanan*) in northern Sweden. The stretch that is considered is between Kiruna (Sweden) and Narvik (Norway) as in Figure 4.

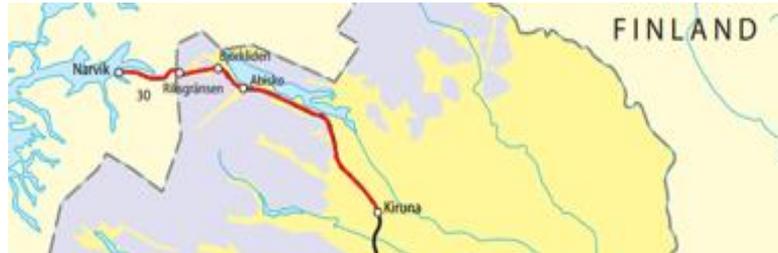
Figure 4 - Map of the Iron Ore line with the stretch between Narvik and Kiruna.

The data was provided by the Swedish National Transport Administration (*Trafikverket*). It consists of the following information:
- Signalling blocks
- Waiting stations
- Travel time between signalling blocks for different scenarios (SF, SS, FS, FF)
- Capacity of waiting stations

In this study, we consider, for the sake of simplicity, only one type of trains having the same speed properties for the different requested train paths.

The data includes a list of 32 train requests from a typical weekday on the line. There are 6 passenger trains are operated by SJ AB and the remaining 26 freight trains are operated by three different freight operators: Green Cargo, Hector Rail and MT AB. The train request information includes:
- Departure and arrival stations
- Ideal departure time
- Latest arrival time at final destination

For each departure, we assume a departure time window of one hour. This means that trains are allowed to depart at the earliest 30 min before the ideal departure time and 30 min at the latest.

Since it was difficult to get data about the valuation of each train request, we used a simplified valuation function. It is a triangular function that has a peak at the ideal departure time and decreases linearly within the departure window. This refers to the definition in Section 2.2 of a utility value $v_p$ associated with each path $p \in \mathcal{P}_r$. The value of the path may depend on how close it is to the ideal departure or arrival time, whether the train must be sided along the line in order to way for other trains thus adding to total travel time. In this particular case a distinction is mainly between freight and passenger services. The peak value is (arbitrarily) set to 500 for freight train requests and 1000 for passenger train requests which are generally more valuable when they depart on time. Figure 5 illustrates the valuation function.

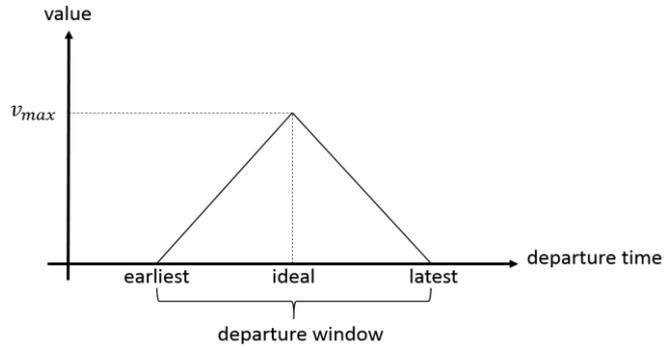
Figure 5 - Simplified valuation function of train requests

In addition, we use a minimal compulsory waiting time of 2min in stations as well as a 3min headway between trains as a blocking rule to ensure a certain safety distance.

In order to check the models on different problem instances, we constructed four different test cases from the given data. Each test case corresponds to the problem of scheduling 32 trains, 26 freight trains with a peak value of 500 and 6 passenger trains with 1000. Table 1 lists the different instances and characteristics.

Table 1 - Test cases and their characteristics

| Test | Termini | # of stations | # of blocks |
|---|---|---|---|
| S1 | Narvik - Bjørnfjell | 5 | 14 |
| S2 | Kiruna - Vassijaure | 7 | 23 |
| S3 | Kiruna - Torneträsk | 14 | 51 |
| S4 | Narvik - Kiruna | 19 | 70 |

The tests cases correspond to an increasingly long stretch of the Iron Ore line. The requested train paths are as previously described and are the same in the different test cases except that the trains depart from and arrive to different termini depending on the test case.

### 4.3 Results and Discussions

The tests were executed on a remote computer with two processors Intel(R) Xeon(R) CPU E5645. Each processor has a clock frequency of 2.40 GHz and 12 MB cache memory; the RAM memory is 80 GB.

The models have several parameters and initialisations to be set before starting the execution. The values that were used for those parameters are given in Table 2.

Table 2 - Algorithm parameters and their values

| Parameter | Value |
|---|---|
| Time discretization step | 30 seconds |
| Step quality threshold | $m_L = 0.1\ (= 10\%)$ |
| Initial step control value | $u_0 = 1$ |
| Minimal step control | $u_{min} = 10^{-10}$ |
| Maximal number of iteration | $k_{max} = 200$ |
| Initial prices | $\mu_0 = 0$ |
| Tolerance (stopping condition) | $\epsilon = 10^{-13}$ |

Both models were tested under the same conditions, i.e. same machine, parameters and input data. Figure 6 shows the comparison between the dual objective in the aggregated and disaggregate approaches for the test cases $S1 - S4$.

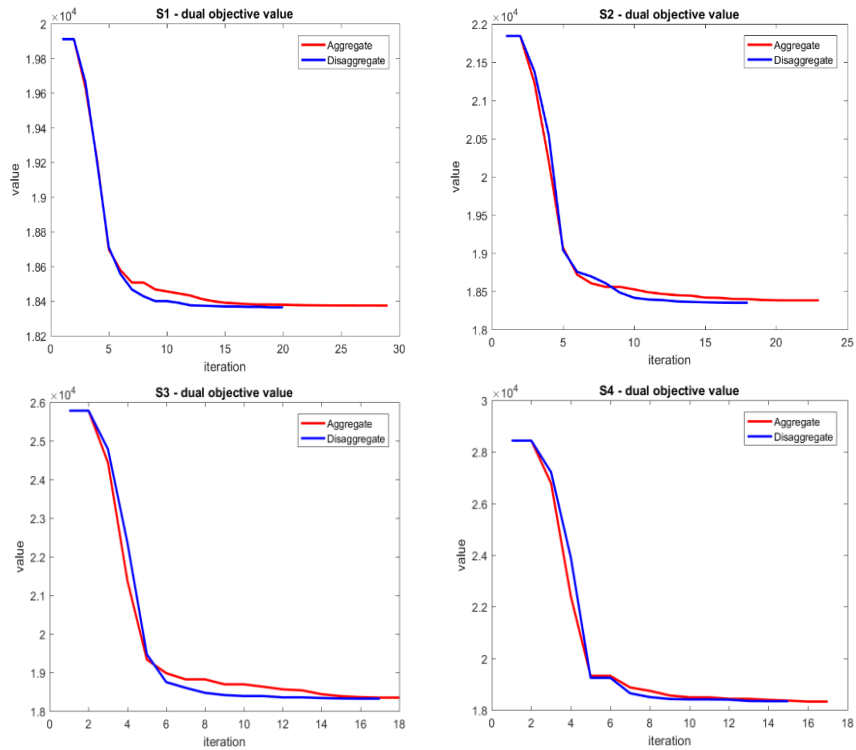

Figure 6- dual objective for the two approaches in the test cases $S1 - S4$ (left to right, up to down)

In the four test cases, the optimization of the dual objective function has a similar behaviour for the two approaches in the first iterations. However, after a certain number of iterations, the minimization in the disaggregate approach becomes faster as more information is collected in the iteration bundle. This leads to a faster convergence using the disaggregate approach.

Each test case leads to an optimization model with a different size, for instance a different number of constraints or variables. The model characteristics of the dual optimization problem (D) are presented in Table 3 for each test case.

Table 3 - Test cases and dual optimisation model characteristics

| Test | # of variables (= # of constraints) |
|------|-------------------------------------|
| S1   | 40 208                              |
| S2   | 66 056                              |
| S3   | 146 472                             |
| S4   | 201 040                             |

Table 4 gives the numerical results regarding the computational time for the test cases. The initialization time corresponds to the time needed to construct the train movement graph. In this step, the algorithm constructs a large space-time graph that includes all the possible train paths. The graph construction time increases when the number of blocks and stations increase or when the time discretization step decreases. This is however done once for each scenario and hence the time is the same for the two approaches. The execution time corresponds to the time the model takes to compute the solution before stopping. It consists of the time needed to compute several shortest path problems over the previously constructed large graph. It also includes the time needed to solve the quadratic optimisation problem in the bundle method for each iteration.

Table 4- Startup and convergence time in the test cases S1 – S4

| Test | Initialization time (in min) | Execution time – aggregate (in min) | Execution time – disaggregate (in min) | Execution time improvement * (in %) |
|------|------------------------------|-------------------------------------|----------------------------------------|-------------------------------------|
| S1   | 26.44                        | 40.36                               | 24.11                                  | 40.3                                |
| S2   | 39.17                        | 49.71                               | 27.13                                  | 45.4                                |
| S3   | 213.19($\approx$ 3.5h)       | 209.09($\approx$ 3.5h)              | 169.72($\approx$ 2.8h)                 | 18.8                                |
| S4   | 374.49($\approx$ 6.2h)       | 334.94($\approx$ 5.6h)              | 253.80($\approx$ 4.2h)                 | 24.2                                |

* The improvement is relative to the convergence time of the aggregate approach.

The computational results in Table 4 show that the disaggregation approach yields an improvement of up to 45% compared to the aggregate in computing the solution to the dual optimisation problem. In addition to this improvement in the execution time, the disaggregation approach is parallelizable, i.e. multiple processors can be used to solve the shortest path subproblems parallelly. This will further improve the computation time and is suggested for future work.

In contrast to Table 4, Figure 6 does not clearly show a faster convergence in the disaggregate approach. More complex problem instances would reveal more clearly the convergence properties of the disaggregate bundle method.

Testing the two approaches provided also additional useful information. For instance, the solution to the dual problem corresponds to an optimal pricing of the infrastructure in space-time which is visualized in
Figure 7 for the first test case S1. Some blocks during some time periods have higher prices than others. These are the blocks corresponding to the requested origin and destination stations during the requested departure and arrival times. The final prices are generally similar between the two approaches with the disaggregate approach yielding slightly higher prices. This slight difference is due to the convergence tolerance that is used in the stopping condition.

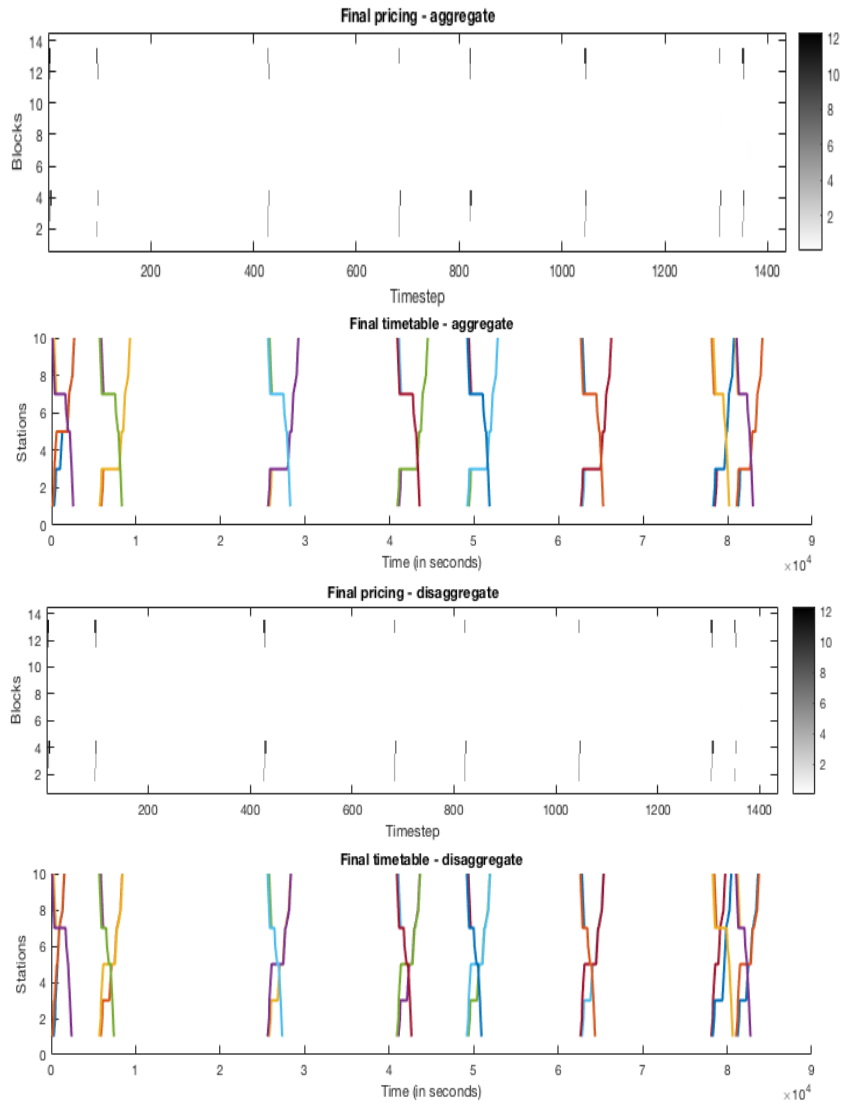

Figure 7- The final pricing and timetable for the two approaches from test case S1

In each iteration, the shortest path algorithm is called to find the best paths for a given pricing. This leads to the generation of a set of potentially optimal paths in each iteration. Therefore, each train request has a set of paths which contains an optimal path (including the null path, i.e. cancelling the train) that will be selected in the final optimal timetable. An example of such as set of path is illustrated in the following figure (from the first request in the test case S4).

The model chooses a new path in each iteration based on the pricing of the infrastructure (Figure 7) by avoiding expensive paths. Figure 8 shows the set of generated paths for one of the train requests (first request in test case S4). It is important to notice that the disaggregate approach yields more paths than the aggregate one which means that it covers more possible solutions in the feasible space. This is very beneficial for a future implementation of a Branch & Bound (B&B) method. These paths will be used later in B&B to find a better feasible solution of the original problem since more potentially optimal paths were generated.

It is also important to note that the null path is always part of such a set which allows to cancel requested trains.

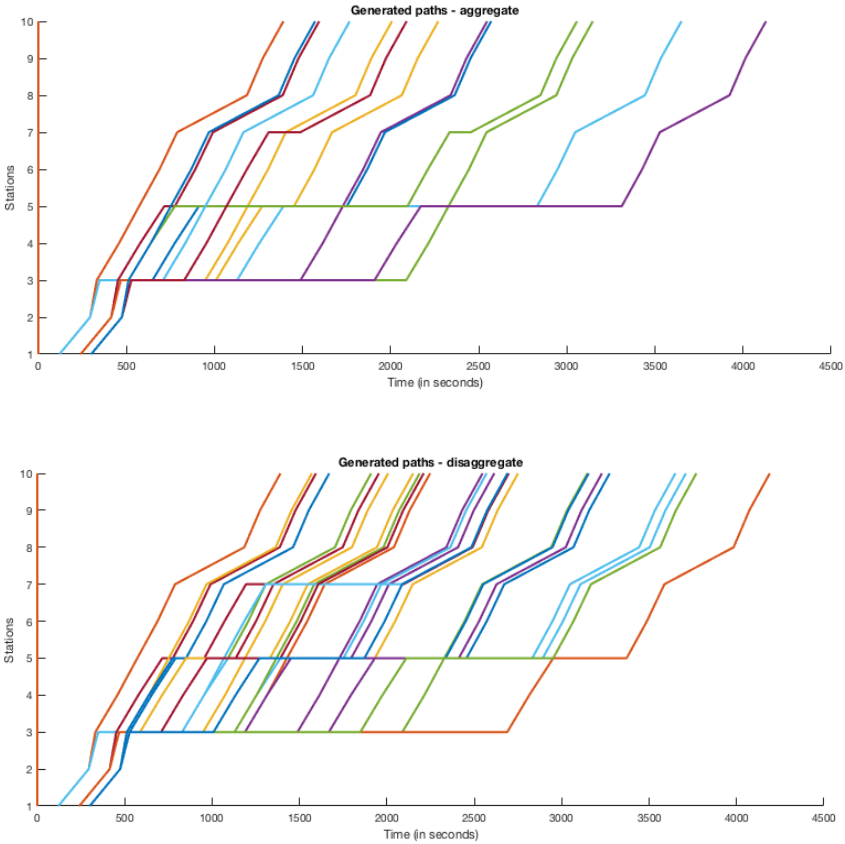

Figure 8- Set of generated paths for the same requested train path using the two approaches (from the first requested train path in test case S1).

## 5  Concluding remarks and Suggestions for Future Work

Lagrangian relaxation is a classical method to solve large scale optimization problem such as the train timetabling problem. The optimization of a dual objective function is often needed in order to gain diverse information about the original problem. Such information can be, for instance, the lower or upper bound of the objective value, Lagrangian multipliers of the relaxed constraints.

Being able to quickly solve the dual problem is very important especially if it is iteratively solved within another optimization model such as Branch & Bound. This paper investigated the potential of a disaggregate approach to speed up the bundle based solution methods for the train timetable problem. This approach is compared with the classical aggregate bundle method. We have demonstrated that the novel approach shows better computational performances in solving the dual problem arising from the relaxation of a train timetabling problem. The novel approach yields faster convergence compared to the classical approach.

In the example given, the value of being allocated a train path has been arbitrarily set at different levels for freight and passenger services. It has been suggested that this input to the optimisation exercise emanates from an explicit bidding process where different operators define the paths they request and the value function of being allocated a path. The utility value $v_p$ associated with the allocated path $p \in \mathcal{P}_r$ then specifies the operator's benefit of being able to run each service. Solving the track allocation problem more generally comprises two components; the optimisation problem which is addressed by the present paper; and the valuation problem which can be handled by operators submitting bids for each path. The latter problem was addressed in Nilsson (2002).

These results have a positive impact on the execution time of train timetabling optimization models that are based on Lagrangian relaxation. To show this, we suggest to investigate the use of the disaggregate approach in a complete Branch & Bound based timetabling model. This has a potential to improve the quality of the final feasible timetable.

## References

BRÄNNLUND, U., LINDBERG, P. O., NÕU, A. & J.-E, N. 1998. Railway Timetabling using Lagrangian Relaxation. *Transportation Science,* 32**,** 358-369.

CAPRARA, A., FISCHETTI, M. & TOTH, P. 2002. Modeling and solving the train timetabling problem. *Operations Research,* 50**,** 851-861.

CAPRARA, A., MONACI, M., TOTH, P. & GUIDA, P. L. 2006. A Lagrangian heuristic algorithm for a real-world train timetabling problem. *Discrete Applied Mathematics,* 154**,** 738-753.

CORMEN, T. H., LEISERSON, C. E. & RIVEST, R. L. 2009. *Introduction to Algorithms (3),* Cambridge, US, The MIT Press.

JAMILI, A., SHAFIA, M. A., SADJADI, S. J. & TAVAKKOLI-MOGHADDAM, R. 2012. Solving a periodic single-track train timetabling problem by an efficient hybrid algorithm. *Engineering Applications of Artificial Intelligence,* 25**,** 793-800.

KIWIEL, K. C. 1990. Proximity control in bundle methods for convex nondifferentiable minimization. *Mathematical Programming,* 46**,** 105-122.
MATHWORKS. 2016. *Introducing MEX Files* [Online]. source MEX fileC, C++, or Fortran source code file. Available: http://se.mathworks.com/help/matlab/matlab_external/introducing-mex-files.html [Accessed October 2015].
NILSSON, J.-E. 2002. Towards a welfare enhancing process to manage railway infrastructure access. *Transportation Research Part A,* 36**,** 419–436.
SCHLECHTE, T. 2012. *Railway Track Allocation: Models and Algorithms.* University of Berlin.
XU, X., LI, K., YANG, L. & YE, J. 2014. Balanced train timetabling on a single-line railway with optimized velocity. *Applied Mathematical Modelling,* 38**,** 894-909.